\documentclass[preprint,review,12pt]{elsarticle}

\usepackage{amssymb}

\journal{journal}
\begin{document}

\begin{frontmatter}
\title{Robin's inequality,Lagarias' criterion, and Riemann hypothesis} 
\author{Ahmad Sabihi \corref{cor1} \cortext[cor1] {Corresponding author}}
\ead{Corresponding author: sabihi2000@yahoo.com,Isfahan city, Iran}
\address{Shahid Ashrafi Esfahani University}
\begin{abstract}
In this paper, we make use of Robin and Lagarias' criteria to prove Riemann hypothesis. The goal is, using Lagarias criterion for $n\geq 1$ since Lagarias criterion states that Riemann hypothesis holds if and only if the inequality $\sum_{d|n}d\leq H_{n}+\exp(H_{n})\log(H_{n})$ holds for all $n\geq 1$. Although, Robin's criterion is used as well. Our approach breaks up the set of the natural numbers into the two main subsets. The first subset is $\{n\in \mathbb{N}| ~ 1\leq n\leq 2(3\times5\dots\times331)^{2}\}$. The second one is $\{n\in \mathbb{N}| ~ n\geq 2(3\times5\dots\times331)^{2}\}$. In our proof, the second subset is decomposed again into the three sub-subsets including odd numbers and the two groups of the even numbers. Then,each group of the even numbers is expressed by an odd integer class number set. Finally, mathematical arguments are stated for each odd integer class number set. Odd integer class number set is introduced in this paper. Since the Lagarias criterion holds for the first subset regarding computer aided computations and Thomas Morrill's paper, we do prove it for the second subset using both Lagarias and Robin's criteria and mathematical arguments. It then follows that Riemann hypothesis holds as well. Essential keys of the proof for large numbers are theorem 1 proving 
$\sigma(m)<\frac{1}{2}e^{\gamma}m \log\log(2m)$ for odd numbers $m\geq (3\times5\dots\times331)^{2}$, lemma 9 and lemma 10 proving $e^{\gamma}(1-\frac{1}{p_{1}})\dots (1-\frac{1}{p_{n}})\log\log(2p_{1}\dots p_{n})<2$ for $n\geq 1$ and $e^{\gamma}(1-\frac{1}{p_{1}})\dots (1-\frac{1}{p_{n}})\log\log(2p^{2}_{1}\dots p^{2}_{n})>2$ for $n\geq 66$.    
\end{abstract}
\begin{keyword}
 Elementary number theory; Analytic number theory; Sum of divisors function; Robin's criterion; Lagarias' criterion; Odd integer class number set\\
\textbf{MSC 2010}:11A05;11A07;11A41;11N05;11M26
\end{keyword}
\end{frontmatter}
\section{Introduction}
Riemann hypothesis (RH) is the most importantly and influentially unsolved problem since Riemann proposed it in a number theory paper in 1859 \cite{BR}. The problem is regarded as one of the most useful and applicable ones not only in pure mathematics, but also in many other scientific fields such as physics and engineering. It's said that the solution of RH would imply immediately the solutions of many other unsolved problems in number theory.\\ 
\indent  The  lemmas and theorems for the proof of RH are presented in Section 2. The proofs are given in Section 3. In Section 2, we break up the entire natural numbers set into the two main subsets.The first subset is $\{n\in \mathbb{N}| ~ 1\leq n\leq 2(3\times5\dots\times331)^{2}\}$. The second one is $\{n\in \mathbb{N}| ~ n\geq 2(3\times5\dots\times331)^{2}\}$.The second subset is decomposed again into the three sub-subsets including odd numbers and the two groups of the even numbers. Then, each group of the even numbers is expressed by an odd integer class number set. Odd integer class number sets are defined in Lemma 6.The goal is, using Lagarias criterion for the proof of RH. Lagarias criterion states that Riemann hypothesis holds if and only if the inequality $\sum_{d|n}d\leq H_{n}+\exp(H_{n})\log(H_{n})$ holds for whole of $n\geq 1$. We make use of Robin's criterion as well.\\ \indent Finally, one takes a conclusion that RH follows from a combination of Lagarias and Robin's criteria for all the natural numbers greater than or equal to 5041. For further research in distribution of primes and Riemann hypothesis refer to references such as \cite{HH},\cite{T},\cite{AW},\cite{PB},\cite{BM},\cite{HD},\cite{NK},\cite{NPK},\cite{KB}.\\\\
\setcounter{equation}{0}
\section{Lemmas, and theorems} \label{sc2}
\textit{\textbf{Lemma 1}. Robin's criterion \cite{R}: States that RH is true if and only if the inequality
\begin{equation}
\sigma(n)< e^{\gamma}n\log\log n, ~~~~~~ for~ n\geq 5041
\end{equation}
holds.}\\
\textit{\textbf{Lemma 2}. Grytczuk's theorem 1 \cite{AG}: Let $n=2m$, $\gcd(2,m)=1$ and $m=\prod_{j=1}^{k} p_{j}^{\alpha_{j}}$. Then for all odd positive integers $m>\frac{3^{9}}{2}$, we have 
\begin{equation}
\sigma(2m)<\frac{39}{40}e^{\gamma}2m\log\log 2m<e^{\gamma}2m\log\log 2m, ~~~~~~ for~ 2m\geq 19686>3^{9}
\end{equation}
and 
\begin{equation}
\sigma(m)<e^{\gamma}m\log\log m, ~~~~~~ for~ m\geq 9843>\frac{3^{9}}{2}
\end{equation}}\\
\textit{\textbf{Lemma 3}. Choie's et.al theorem 1.2 \cite{YC}: Any odd positive integer $n$ distinct from 1,3,5 and 9 is in ${\cal R}$, where ${\cal R}$ denotes the set of integers $n\geq 1$ satisfying 
\begin{equation}
\sigma(n)<e^{\gamma}n\log\log n
\end{equation}}
\textit{\textbf{Lemma 4}. Lagarias' theorem 1.1 \cite{L2}: Let $H_{n}=\sum_{j=1}^{n}\frac{1}{j}$. Show that, for each $n\geq 1$, we have
\begin{equation}
\sum_{d|n}d\leq H_{n}+\exp(H_{n})\log(H_{n})
\end{equation} 
with equality only for $n=1$, where $H_{n}$ denotes $nth$ harmonic number. Then correctness of the inequality (2.5) is equivalent to correctness of RH.}\\
\textit{\textbf{Lemma 5}. Lagarias' lemma 3.1 \cite{L2}: For all $n\geq 3$, we have
\begin{equation}
 e^{\gamma}n\log\log n< H_{n}+\exp(H_{n})\log(H_{n})
\end{equation}}
\textit{\textbf{Lemma 6}. 
The set of all the even numbers greater than or equal to $N$ is equivalent to the union of the odd integer class number sets of 1 to infinity. This means that
\begin{equation}
\left\{n\geq N | n~is~an~even~number\right\}=\cup_{m=1}^{\infty}CL(m)
\end{equation}}
\textit{\textbf{Definition 1}. Odd integer class number set. Odd integer class number set is a class of even integers $n=2^{\alpha}m\geq N$ for an arbitrarily even number $N$, where $m$ denotes an odd number demonstrating an odd class number and $\alpha$ denotes a natural number. The odd integer class number set of $m$ is defined as the set
\begin{equation}
CL(m)=\left\{2^{\alpha}m\geq N |  m~is~ a~constantly~ odd~ number~ and~\alpha\geq\alpha_{0} \right\}
\end{equation}
 and $\alpha_{0}$ denotes the smallest positive integer to hold $2^{\alpha}m\geq N$. For example: Odd integer class number set of 3 is the set (provided that $2^{\alpha_{0}}\times3\geq N$, where $N$ is also an even number).
\begin{equation}
CL(3)=\left\{2^{\alpha_{0}}\times3, 2^{\alpha_{0}+1}\times3, 2^{\alpha_{0}+2}\times3\dots \right\}
\end{equation}}
\textit{\textbf{Lemma 7} (Rosser et. al \cite{SBR})\\
Let $p_{n}>285$, then
\begin{eqnarray}
\frac{e^{-\gamma}}{\log p_{n}}(1-\frac{1}{2\log^{2}p_{n}})<\prod_{p\leq p_{n}}(1-\frac{1}{p})<\frac{e^{-\gamma}}{\log p_{n}}(1+\frac{1}{2\log^{2}p_{n}})
\end{eqnarray}}
Note that $p_{1}=3$ is considered the origin of counting the odd primes in this lemma.\\
\textit{\textbf{Lemma 8} (Rosser et. al \cite{SBR})\\
Let $p_{n}\geq 41$, then
\begin{eqnarray}
p_{n}(1-\frac{1}{\log p_{n}})<\log(2 p_{1} p_{2}\dots p_{n})<p_{n}(1+\frac{1}{2\log p_{n}})
\end{eqnarray}}
Note that $p_{1}=3$ is considered the origin of counting the odd primes in this lemma.\\
\textit{\textbf{Proposition 1}\\
Let $g(x)\geq 0$, $g'(x)>0$,$f(x)$ and $f'(x)$ be defined and differentiable on $f,g: [0,+\infty[\longrightarrow \mathbb{R}$. None of the functions $f(x)$ and $g(x)$ are sinusoidal or oscillating or having changed signs. Let $f(x)=O(g(x))$, be defined on the entire domain $f,g: [0,+\infty[ \longrightarrow \mathbb{R}$ and $g(x_{0})=0$, where $x_{0}$ belongs to $[0,+\infty[$  or be defined for sufficiently large values $x$ and $\lim_{x\rightarrow +\infty}g(x)=+\infty$ for both cases, then
\begin{center}
$f'(x)=O(g'(x))$
\end{center}
where $O$ denotes big O-notation function.}\\\\
\textit{\textbf{Lemma 9}\\
Let $n\geq 1$, then the function
\begin{eqnarray}
RO_{1}(n)=e^{\gamma}(1-\frac{1}{p_{1}})\dots (1-\frac{1}{p_{n}})\log\log(2p_{1}\dots p_{n})
\end{eqnarray} 
is a strictly increasing one for consecutive primes $p_{1}$ to $p_{n}$.}
Note that $p_{1}=3$ is considered the origin of counting the odd primes in this lemma.\\ 
\textit{\textbf{Lemma 10}\\
Let $n\geq 66$, then for consecutive primes $p_{1}$ to $p_{n}$, we have 
\begin{eqnarray}
RO_{2}(n)=e^{\gamma}(1-\frac{1}{p_{1}})\dots (1-\frac{1}{p_{n}})\log\log(2p^{2}_{1}\dots p^{2}_{n})>2
\end{eqnarray} 
and 
\begin{eqnarray}
RO_{1}(n)<2
\end{eqnarray}}
Note that $p_{1}=3$ is considered the origin of counting the odd primes in this lemma.\\
\textit{\textbf{Lemma 11}. RH is true if and only if the following inequalities for $n\geq 5041$ hold \cite{L2}
\begin{equation}
\sigma(n)< e^{\gamma}n\log\log n< H_{n}+\exp(H_{n})\log(H_{n})
\end{equation}}
\textit{\textbf{Theorem 1}. Let 
\begin{equation}
\sigma(m)<\frac{1}{2}e^{\gamma}m \log\log(2m)
\end{equation} 
be true for $m=m_{0}=(3\times5\times 7\times\dots \times331)^{2}$, then it is correct for all the odd numbers $m\geq m_{0}$}\\
\textit{\textbf{Theorem 2}. Lagarias criterion holds for $\{n\in \mathbb{N}| ~ 1\leq n\leq 2(3\times5\dots\times331)^{2}\}$}.\\
\textit{\textbf{Theorem 3}. Lagarias and Robin's criteria hold for the odd numbers $\{n\in \mathbb{N}| ~ n> 2(3\times5\dots\times331)^{2}\}$}.\\  
\textit{\textbf{Theorem 4}. Lagarias and Robin's criteria hold for all the odd integer class number sets of $m$, where $m$ is an odd integer $1\leq m<(3\times5\dots\times331)^{2}$ and all the even numbers $n\geq 2(3\times5\dots\times331)^{2}$.}\\
\textit{\textbf{Theorem 5}. Lagarias and Robin's criteria hold for all the odd integer class number sets of $m$, where $m$ is an odd integer $m\geq (3\times5\dots\times331)^{2}$ and all the even numbers $n\geq 2(3\times5\dots\times331)^{2}$.}\\
\setcounter{equation}{0}   
\section{Proofs of the lemmas and theorems} \label{sc4}
\subsection{\textbf{Proof of Lemma 1}} The proof is found in the Robin's paper \cite{R}. 
\subsection{\textbf{Proof of Lemma 2}} The proof is found in the Grtytczuk's paper \cite{AG}.
\subsection{\textbf{Proof of Lemma 3}} The proof is found in the Choie's et. al paper \cite{YC}.
\subsection{\textbf{Proof of Lemma 4}} The proof is found in the Lagarias' paper \cite{L2}.
\subsection{\textbf{Proof of Lemma 5}} The proof is found in the Lagarias' paper \cite{L2}.
\subsection{\textbf{Proof of Lemma 6}} Trivially, the proof is very easy. Regarding definition1, let an even number as $n$, where $\alpha_{0}$ denotes the smallest positively integer value that $2^{\alpha_{0}}m\geq N$ holds, then for every $\alpha\geq \alpha_{0}$, we have  
$2^{\alpha}m\geq 2^{\alpha_{0}}m\geq N$. This states the odd integer class number set $CL(m)$ corresponding to $m$. If we establish the set $CL(m)$ for each $m=1,3,5,...$, and make the set of their union, then we have $\cup_{m=1}^{\infty}CL(m)$. On the other hand, let $m_{0}$ and $m_{1}$ be two distinctly odd class numbers and their corresponding even numbers be $2^{\alpha_{0}}m_{0}$ and $2^{\alpha_{1}}m_{1}$, respectively. There is no equality $2^{\alpha_{0}}m_{0}=  2^{\alpha_{1}}m_{1}$ since equivalency contradicts being odd $m_{0}$ or $m_{1}$. Therefore, $\cup_{m=1}^{\infty}CL(m)$ must cover and represent all the even numbers greater than or equal to $N$.
\subsection{\textbf{Proof of Lemma 7}}
The proof is made in Rosser and Schoenfeld's paper \cite{SBR}.
\subsection{\textbf{Proof of Lemma 8}}
The proof is made in Rosser and Schoenfeld's paper \cite{SBR}.
\subsection{\textbf{Proof of Proposition 1}}
If $f(x)=O(g(x))$, would imply that  
\begin{center}
$\lim_{x \rightarrow \infty}sup \frac{|f(x)|}{g(x)}<\infty$
\end{center}
or
\begin{center}
$|f(x)|<Ag(x)$
\end{center}
where $A>0$ and denotes a constant value independent of $x$.\\
The problem is proved for the two cases: 1-for when $|f(x)|<Ag(x)$ holds on the entire domain $f,g$. 2-for when $|f(x)|<Ag(x)$ holds only for sufficiently large values $x$.\\  
\textbf{Case 1:}\\
Just, we prove that 
\begin{center}
$|f'(x)|<Ag'(x)$
\end{center}
holds on the domain $f,g: [0,+\infty[\longrightarrow \mathbb{R}$. If $f(x)$ is a constant function, then $|f'(x)|=0<Ag'(x)$. On the other hand, function $f(x)$ is assumed be unbounded one and the case of a bounded function is given to be followed in lately proof. Therefore, we ignore these in both cases 1 and 2.
The proof is made by \textit{reductio ad absurdum}.\\
If the inequality $|f'(x)|<Ag'(x)$ is assumed to be false on the entire domain $f,g$, then we have by contradiction (if it is not false in some of points, then we ignore them because the inequality $|f'(x)|<Ag'(x)$ holds in such points)
\begin{center}
$|f'(x)|\geq Ag'(x)$
\end{center}
holds on all the domain $f,g$.\\ 
If both $f'(x)>0$ and $f(x)>0$, we have
\begin{center}
$f'(x)\geq Ag'(x)$
\end{center}
and are able to integrate the both hand-sides over the entire domain $(x,x_{0}\in D_{f,g})$ and find
\begin{center}
$\int_{x_{0}}^{x}f'(u)du\geq \int_{x_{0}}^{x}Ag'(u)du$
\end{center}
and
\begin{center}
$(f(x)-f(x_{0}))\geq A(g(x)-g(x_{0}))$
\end{center}
since $g(x_{0}=0)=0$, then the inequality $|f(x)|<Ag(x)$ implies that $f(x_{0}=0)=0$ as well and we have 
\begin{center}
$f(x)\geq Ag(x)$
\end{center}
on the domain $f,g$ and leads a contradiction with assumption since we must have $f(x)\leq Ag(x)$ and the proof is completed.\\
If $f'(x)<0$ and $f(x)>0$, then by contradiction, we have
\begin{center}
$f'(x)\leq -Ag'(x)$
\end{center} 
Similarly to the above, integrating it yields $(x,x_{0}\in D_{f,g})$
\begin{center}
$\int_{x_{0}}^{x}f'(u)du\leq -\int_{x_{0}}^{x}Ag'(u)du$
\end{center}
and we have with $x_{0}=0$
\begin{center}
$f(x)\leq -Ag(x)$
\end{center}
on the domain $f,g$ and leads a contradiction with assumption since we must have $f(x)\leq Ag(x)$  and the proof is completed. Similarly to the case $f(x)>0$, the proof is valid when $f(x)<0$ as well.\\
\textbf{Case 2:}\\
Let $|f'(x)|<Ag'(x)$ be a false statement for $x>x_{N}$, where $x_{N}$ is a sufficiently large number, then regarding contradiction, we should have $|f'(x)|\geq Ag'(x)$ for $x>x_{N}$.\\
If both $f(x)$ and $f'(x)>0$,we should have
\begin{center}
$f'(x)\geq Ag'(x)$
\end{center}
 for $x>x_{N}$. If we integrate both hand-sides of the above inequality between $x$ and  $x_{N}$, we find the following inequality
\begin{center}
$h(x)=f(x)-Ag(x)\geq h(x_{N})=f(x_{N})-Ag(x_{N})$
\end{center}
On the other hand, $h(x)$ is an unbounded function since regarding the proposition's assumption $g(x)$ tends to infinity when $x$ so does. This makes the inequality $h(x)\geq h(x_{N})$ incorrect since $h(x)$ tends to $-\infty$, but $h(x_{N})$ is a bounded value and fixed. This contradicts our problem's assumption. Thus, $|f'(x)|<Ag'(x)$ is correct. \\
If $f(x)>0$ and $f'(x)<0$, then we should have $f'(x)\leq -Ag'(x)$. After integrating both hand sides of the inequality between $x$ and  $x_{N}$, we have
\begin{center}
$f(x)\leq f(x_{N})+Ag(x_{N})-Ag(x)$
\end{center}   
and regarding the problem's assumption, we have
\begin{center}
$f(x)<Ag(x)$
\end{center}   
Summing out the above inequalities together,yields
\begin{center}
$2f(x)\leq f(x_{N})+Ag(x_{N})$
\end{center}  
where restricts $f(x)$, since $f(x)$ is not a bounded function (the cases that $f(x)$ is a fixed or a bounded one are ignored as mentioned at the beginning of the proof. Being bounded $f(x)$ contradicts the inequality $f(x)\leq f(x_{N})+Ag(x_{N})-Ag(x)$ since here $f(x)>0$, but the right hand of the inequality goes to the negative state and gets contradicted), but $f(x_{N})+Ag(x_{N})$ is a bounded value. Again, we find out that $|f'(x)|<Ag'(x)$ is correct.\\
If both $f(x)$ and $f'(x)<0$, then the inequality $|f(x)|<Ag(x)$ implies that $f(x)>-Ag(x)$. This means that holding it with the inequality $f(x)\leq f(x_{N})+Ag(x_{N})-Ag(x)$ simultaneously, leads to a contradiction since we find
\begin{center}
$0<f(x)+Ag(x)<f(x_{N})+Ag(x_{N})$
\end{center}
because $f(x)+Ag(x)$ is an unbounded function when $x$ tends to infinity but $f(x_{N})+Ag(x_{N})$ is a bounded value. All of the above reasons yields $f'(x)\leq -Ag'(x)$. Therefore, $f'(x)=O(g'(x))$ holds.
\subsection{\textbf{Proof of Lemma 9}}
We directly investigate that the lemma is correct for all $p_{n}< e^{130000}$. The proof for $x\geq p_{n}> e^{130000}$ is made by analytic reasoning as follows:\\
Just, we compute $\log \theta(x)=\log\log_{p_{n}\leq x}(2p_{1}\dots p_{n})$ regarding the relation (2.12) and using Abel Summation Formula (\cite{WE}, page 10, theorem 1.6).\\
Let $a_{n}=1$ be a constant sequence,$A(x)=\sum_{p\leq x}a_{n}=\pi(x)$ and $\phi_{1}(p)=\log p$. Abel summation formula asserts 
\begin{eqnarray}
\theta(x)=\sum_{p\leq x}\log p=\pi(x)\log x-\int_{2}^{x}\frac{\pi(u)}{u}du
\end{eqnarray}
Thus,
\begin{eqnarray}
\log\theta(x)=\log\log_{p_{n}\leq x}(2p_{1}\dots p_{n})=\log\left(\pi(x)\log x-\int_{2}^{x}\frac{\pi(u)}{u}du\right)
\end{eqnarray}
First of all, to compute the expression $B=\prod_{p\leq p_{n}\leq x}(1-\frac{1}{p})$ of (2.12), we have to let $p_{1}=2$ and compute (3.3) assuming it 
\begin{eqnarray}
\log B=\sum_{p\leq p_{n}\leq x}\log(1-\frac{1}{p})
\end{eqnarray}
where $\phi_{2}(p)=\log(1-\frac{1}{p})$, and again ,we have $A(x)=\sum_{p\leq x}a_{n}=\sum_{p\leq x}1=\pi(x)$ and regarding Abel summation formula yields 
\begin{eqnarray}
\log B=\pi(x)\log(1-\frac{1}{x})-\int_{2}^{x}\frac{\pi(u)}{u(u-1)}du
\end{eqnarray}
that implies that
\begin{eqnarray}
B=e^{\left\{\pi(x)\log(1-\frac{1}{x})-\int_{2}^{x}\frac{\pi(u)}{u(u-1)}du\right\}}=(1-\frac{1}{x})^{\pi(x)}e^{-\int_{2}^{x}\frac{\pi(u)}{u(u-1)}du}
\end{eqnarray}
Let us change the argument $n$ in $RO_{1}(n)$ given in (2.12) into $x$ as follows letting $p_{1}=3$:
\begin{eqnarray}
C(x)=e^{\gamma}(1-\frac{1}{p_{1}})\dots (1-\frac{1}{p_{n}})_{2<p_{n}\leq x}\log\log(2p_{1}\dots p_{n})_{2<p_{n}\leq x}=\nonumber\\ e^{\gamma}\prod_{2<p\leq p_{n}\leq x}(1-\frac{1}{p})\log \theta(x)=2e^{\gamma}(1-\frac{1}{x})^{\pi(x)}_{2<p\leq x}e^{-\int_{2}^{x}\frac{\pi(u)}{u(u-1)}du}\log \theta(x)
\end{eqnarray}
Since we know the function $(1-\frac{1}{x})^{\pi(x)}_{2<p\leq x}$ (the first term of the right -hand side of (3.6)) is a strictly increasing one, it is enough to prove the following function (the right terms of the right hand side of (3.6)) is a strictly increasing one:
\begin{eqnarray}
Y=\frac{\log \theta(x)}{e^{\int_{2}^{x}\frac{\pi(u)}{u(u-1)}du}}
\end{eqnarray} 
Taking the first derivative of (3.7) and considering only its numerator, regarding (3.2), gives us $Z$. ( Note that the denominator is a positive value and no needs to take it into account)
\begin{eqnarray}
Z=x(x-1)\pi'(x)\log x-\pi(x)\left\{\pi(x)\log x-\int_{2}^{x}\frac{\pi(u)}{u}du\right\}\times\nonumber\\ \left\{\log\left(\pi(x)\log x-\int_{2}^{x}\frac{\pi(u)}{u}du\right)\right\}~~~~~~~~~~~~~~~~
\end{eqnarray} 
Just, we should prove that $Z>0$ for $x\geq e^{130000}$ since we can suppose that $x$ is a sufficiently large number.\\
Prime number theorem asserts \cite{WE}, page 137. 
\begin{eqnarray}
\pi(x)=\int_{2}^{x}\frac{du}{\log u}+O\left(\frac{x}{e^{a\sqrt{\log x}}}\right)
\end{eqnarray} 
and differentiating $\pi(x)$ regarding proposition 1, yields
\begin{eqnarray}
\pi'(x)=\frac{1}{\log x}+O\left(\frac{(1-\frac{a}{2\sqrt{\log x}})}{e^{a\sqrt{\log x}}}\right)
\end{eqnarray}
since we have the case 2 of proposition 1, then $f(x)=O\left(\frac{x}{e^{a\sqrt{\log x}}}\right)$, $f'(x)=O\left(\frac{(1-\frac{a}{2\sqrt{\log x}})}{e^{a\sqrt{\log x}}}\right)$ and $g(x)=\frac{x}{e^{a\sqrt{\log x}}}$, $g'(x)=\frac{(1-\frac{a}{2\sqrt{\log x}})}{e^{a\sqrt{\log x}}}$, $g(x)>0$, $g'(x)>0$, and $\lim_{x \rightarrow +\infty}g(x)=+\infty$.
Substituting (3.9) and (3.10) for (3.8), we have
\begin{eqnarray}
Z=x(x-1)+O\left(x(x-1)\frac{(\log x-\frac{a\sqrt{\log x}}{2})}{e^{a\sqrt{\log x}}}\right)-\nonumber\\ \left\{\int_{2}^{x}\frac{du}{\log u}+O\left(\frac{x}{e^{a\sqrt{\log x}}}\right)\right\}\times~~~~~~~~~~~~~ \nonumber\\ \left\{\log x \int_{2}^{x}\frac{du}{\log u}-\int_{2}^{x}\frac{\pi(u)du}{u}+O\left(\frac{x\log x}{e^{a\sqrt{\log x}}}\right)\right\}\times\nonumber\\ \log\left\{\log x \int_{2}^{x}\frac{du}{\log u}-\int_{2}^{x}\frac{\pi(u)du}{u}+O\left(\frac{x\log x}{e^{a\sqrt{\log x}}}\right)\right\} 
\end{eqnarray}  
Denoting the third term of the (3.11) as $S$ and expanding it,  we have  
\begin{eqnarray}
S=\left\{\int_{2}^{x}\frac{du}{\log u}+O\left(\frac{x}{e^{a\sqrt{\log x}}}\right)\right\}\times~~~~~~~~~~~~~ \nonumber\\ \left\{\log x \int_{2}^{x}\frac{du}{\log u}-\int_{2}^{x}\frac{\pi(u)du}{u}+O\left(\frac{x\log x}{e^{a\sqrt{\log x}}}\right)\right\}\times\nonumber\\ \log\left\{\log x \int_{2}^{x}\frac{du}{\log u}-\int_{2}^{x}\frac{\pi(u)du}{u}+O\left(\frac{x\log x}{e^{a\sqrt{\log x}}}\right)\right\}=\nonumber\\ \left\{(\int_{2}^{x}\frac{du}{\log u})^{2}\log x-(\int_{2}^{x}\frac{du}{\log u})(\int_{2}^{x}\frac{\pi(u)du}{u})\right\}\times\nonumber\\ \log\left\{\log x \int_{2}^{x}\frac{du}{\log u}-\int_{2}^{x}\frac{\pi(u)du}{u}+O\left(\frac{x\log x}{e^{a\sqrt{\log x}}}\right)\right\}+\nonumber\\ \{O\left(\frac{x\log x\int_{2}^{x}\frac{du}{\log u}}{e^{a\sqrt{\log x}}}\right)+O\left(\frac{x\log x\int_{2}^{x}\frac{du}{\log u}}{e^{a\sqrt{\log x}}}\right)+O\left(\frac{x\int_{2}^{x}\frac{\pi(u)du}{u}}{e^{a\sqrt{\log x}}}\right)+\nonumber\\ O\left(\frac{x^{2}\log x}{e^{2a\sqrt{\log x}}}\right)\}\log\left\{\log x \int_{2}^{x}\frac{du}{\log u}-\int_{2}^{x}\frac{\pi(u)du}{u}+O\left(\frac{x\log x}{e^{a\sqrt{\log x}}}\right)\right\}
\end{eqnarray} 
Since 
\begin{eqnarray}
O\left(\frac{x\log x\int_{2}^{x}\frac{du}{\log u}}{e^{a\sqrt{\log x}}}\right)+O\left(\frac{x\log x\int_{2}^{x}\frac{du}{\log u}}{e^{a\sqrt{\log x}}}\right)+O\left(\frac{x\int_{2}^{x}\frac{\pi(u)du}{u}}{e^{a\sqrt{\log x}}}\right)+\nonumber\\ O\left(\frac{x^{2}\log x}{e^{2a\sqrt{\log x}}}\right)=O\left(\frac{x\log x\int_{2}^{x}\frac{du}{\log u}}{e^{a\sqrt{\log x}}}\right)~~~~~~~~~~~~~~~
\end{eqnarray}
because $\frac{x\log x\int_{2}^{x}\frac{du}{\log u}}{e^{a\sqrt{\log x}}}$ is absolutely greater than those of $\frac{x\log x\int_{2}^{x}\frac{du}{\log u}}{e^{a\sqrt{\log x}}}, \frac{x\int_{2}^{x}\frac{\pi(u)du}{u}}{e^{a\sqrt{\log x}}}$ and $\frac{x^{2}\log x}{e^{2a\sqrt{\log x}}}$.
Therefore, regarding (3.12) and (3.13), we have
\begin{eqnarray}
S=\left\{(\int_{2}^{x}\frac{du}{\log u})^{2}\log x-\int_{2}^{x}\frac{du}{\log u}\int_{2}^{x}\frac{\pi(u)du}{u}+O\left(\frac{x\log x\int_{2}^{x}\frac{du}{\log u}}{e^{a\sqrt{\log x}}}\right)\right\}\times\nonumber\\ \log\left\{\log x \int_{2}^{x}\frac{du}{\log u}-\int_{2}^{x}\frac{\pi(u)du}{u}+O\left(\frac{x\log x}{e^{a\sqrt{\log x}}}\right)\right\}~~~~~~~~~
\end{eqnarray}
If we substitute the relation (3.15) based on the Prime Number Theorem for the last term of the right-hand side of (3.14), i.e. to make simple our calculations, we make use of\\ $\theta (x)=x+O\left(\frac{x}{e^{a\sqrt{\log x}}}\right)$ instead of $\left\{\log x \int_{2}^{x}\frac{du}{\log u}-\int_{2}^{x}\frac{\pi(u)du}{u}+O\left(\frac{x\log x}{e^{a\sqrt{\log x}}}\right)\right\}$ in the last term of the right hand side of (3.14).
\begin{eqnarray}
\log\theta(x)=\log\left(x+O\left(\frac{x}{e^{a\sqrt{\log x}}}\right)\right)=\log x+\log\left(1+O(\frac{1}{e^{a\sqrt{\log x}}})\right)
\end{eqnarray} 
then,
\begin{eqnarray}
\log\theta(x)=\log x+\epsilon
\end{eqnarray}
where $\epsilon=\lim_{x\rightarrow \infty}\log\left(1+O(\frac{1}{e^{a\sqrt{\log x}}})\right)$ and tends to zero.\\
(3.14) implies that
\begin{eqnarray}
S=\left\{(\int_{2}^{x}\frac{du}{\log u})^{2}\log x-\int_{2}^{x}\frac{du}{\log u}\int_{2}^{x}\frac{\pi(u)du}{u}+O\left(\frac{x\log x\int_{2}^{x}\frac{du}{\log u}}{e^{a\sqrt{\log x}}}\right)\right\}\times\nonumber\\ \left\{\log x+\epsilon\right\}=
(\int_{2}^{x}\frac{du}{\log u})^{2}(\log x)^{2}-\left\{\int_{2}^{x}\frac{du}{\log u}\int_{2}^{x}\frac{\pi(u)du}{u}\right\}\log x+\nonumber\\ O\left(\frac{x(\log x)^{2}\int_{2}^{x}\frac{du}{\log u}}{e^{a\sqrt{\log x}}}\right)+\epsilon(\int_{2}^{x}\frac{du}{\log u})^{2}\log x-\epsilon\int_{2}^{x}\frac{du}{\log u}\int_{2}^{x}\frac{\pi(u)du}{u}~~~~~
\end{eqnarray}
since $\frac{x(\log x)^{2}\int_{2}^{x}\frac{du}{\log u}}{e^{a\sqrt{\log x}}}$  is greater than $\frac{ x\log x\int_{2}^{x}\frac{du}{\log u}}{e^{a\sqrt{\log x}}}$, we consider only $O\left(\frac{x(\log x)^{2}\int_{2}^{x}\frac{du}{\log u}}{e^{a\sqrt{\log x}}}\right)$ in (3.17).\\
Just, we have
\begin{eqnarray}
Z=x(x-1)+O\left(x(x-1)\frac{(\log x-\frac{a\sqrt{\log x}}{2})}{e^{a\sqrt{\log x}}}\right)-S=\nonumber\\x(x-1)-(\int_{2}^{x}\frac{du}{\log u})^{2}(\log x)^{2}+\left\{\int_{2}^{x}\frac{du}{\log u}\int_{2}^{x}\frac{\pi(u)du}{u}\right\}\log x-\nonumber\\\epsilon(\int_{2}^{x}\frac{du}{\log u})^{2}\log x+\epsilon\int_{2}^{x}\frac{du}{\log u}\int_{2}^{x}\frac{\pi(u)du}{u}+\nonumber\\O\left(\frac{x(\log x)^{2}\int_{2}^{x}\frac{du}{\log u}}{e^{a\sqrt{\log x}}}\right)~~~~~~~~~~~~~~~~~~~~~
\end{eqnarray}
since we know 
\begin{eqnarray}
O\left(\frac{x(\log x)^{2}\int_{2}^{x}\frac{du}{\log u}}{e^{a\sqrt{\log x}}}\right)+O\left(x(x-1)\frac{(\log x-\frac{a\sqrt{\log x}}{2})}{e^{a\sqrt{\log x}}}\right)=\nonumber\\ O\left(\frac{x(\log x)^{2}\int_{2}^{x}\frac{du}{\log u}}{e^{a\sqrt{\log x}}}\right)~~~~~~~~~~~~~~~~~~~
\end{eqnarray}
due to being greater the first term $O$-notation than the second one.\\
To evaluate $\int_{2}^{x}\frac{\pi(u) du}{u}$ given in (3.18) and for sufficiently large number $x>N$ (N is a sufficiently large number as  $e^{130000}$), we have the following identity when x tends to infinity since we know $\frac{x-2}{\log x}<li(x)<\frac{x}{\log x-2}$ for $x\geq e^{4}$ (see lemmas 7,8 of \cite{BROS}) because as $x$ gets sufficiently large and dividing both sides of the inequality by $\frac{x-2}{\log x}$, $li(x)$ tends to $\frac{(x-2)}{\log x}$: 
\begin{eqnarray} 
\int_{2}^{x}\frac{du}{\log u}=\frac{(x-2)}{\log x}+\beta
\end{eqnarray}
where $\beta>0$ and we can also make use of the mean value theorem in integrals.\\
Regarding equivalency $\pi(x)=\frac{x}{\log x}+\alpha$ for sufficiently large numbers, where $\alpha >0$ since we know $\frac{x}{\log x+2}<\pi(x)<\frac{x}{\log x-4}$ for $x\geq 55$ (see theorems 29 and 30,part $A$ of \cite{BROS} and corollary 1 of  \cite{SBR}) for sufficiently large number $x>N$ and dividing both sides of the inequality by $\frac{x}{\log x}$ we have $\pi(x)=\frac{x}{\log x}+\alpha$ then putting (3.20) and $\pi(x)=\frac{x}{\log x}+\alpha$ into (3.21) and regarding corollary 1 of  \cite{SBR}, we have
\begin{eqnarray} 
\int_{2}^{x}\frac{\pi(u) du}{u}=\int_{2}^{17}\frac{\pi(u) du}{u}+\int_{17}^{x}\frac{\pi(u) du}{u}=\frac{(x-2)}{\log x}+\alpha\log \frac{x}{2}+\beta~~~~
\end{eqnarray}
where $\lim_{x \rightarrow \infty}\beta=0$ and $\lim_{x \rightarrow \infty}\alpha=0$. $\alpha$ and $\beta$ are independent of $x$, but tend to zero as $x$ gets large. Therefore, we can omit $\int_{2}^{17}\frac{\pi(u) du}{u}$ versus $\int_{17}^{x}\frac{\pi(u) du}{u}$ and apply corollary 1. This means that $\int_{2}^{x}\frac{\pi(u) du}{u}\sim \int_{17}^{x}\frac{\pi(u) du}{u}$.\\
Just, we substitute (3.20) and (3.21) for (3.18) to investigate if $Z$ is positive or negative as follows: 
\begin{eqnarray} 
Z=x^{2}-x-(\frac{(x-2)}{\log x})^{2}{\log x}^{2}-2\beta\frac{(x-2)}{\log x}(\log x)^{2}-\nonumber\\ \beta^{2}(\log x)^{2}+(\frac{(x-2)}{\log x})^{2}{\log x}+\alpha\log\frac{x}{2}(x-2)+\beta(x-2)+\nonumber\\ \beta(x-2)+\alpha\beta\log\frac{x}{2}\log x+\beta^{2}\log x-\epsilon\frac{(x-2)^{2}}{\log x}-\nonumber\\ 2\epsilon\beta(x-2)-\epsilon\beta^{2}\log x+\epsilon\frac{(x-2)^{2}}{(\log x)^{2}}+\epsilon\alpha\log\frac{x}{2}\frac{(x-2)}{\log x}+\nonumber\\ \epsilon\beta\frac{(x-2)}{\log x}+\epsilon\beta\frac{(x-2)}{\log x}+\alpha\epsilon\beta\log\frac{x}{2}+~~~~~~~~~~~\nonumber\\ \epsilon \beta^{2}-A\beta\frac{x(\log x)^{2}}{e^{a\sqrt{\log x}}}-A\frac{x(x-2)\log x}{e^{a\sqrt{\log x}}}~~~~~~~~~~~~~
\end{eqnarray}
where $-A\beta\frac{x(\log x)^{2}}{e^{a\sqrt{\log x}}}$ and $-A\frac{x(x-2)\log x}{e^{a\sqrt{\log x}}}$ denote the least negative values of big $O$-notation in the last term of (3.18) after substituting (3.20) into big $O$-notation and decomposing it into the two aforementioned terms since $\ O\left(\frac{x(\log x)^{2}\int_{2}^{x}\frac{du}{\log u}}{e^{a\sqrt{\log x}}}\right)\leq A\frac{x(\log x)^{2}\int_{2}^{x}\frac{du}{\log u}}{e^{a\sqrt{\log x}}}$, where $A$ denotes a constant value. Therefore, we can omit the smaller term versus greater one and write as follows:
\begin{eqnarray} 
Z=3x-4+\frac{(x-2)^{2}}{\log x}+\alpha\log\frac{x}{2}(x-2)+2\beta(x-2)+\nonumber\\ \alpha\beta\log\frac{x}{2}\log x+\beta^{2}\log x-2\beta(x-2)\log x-\beta^{2}(\log x)^{2}+\nonumber\\ \epsilon\frac{(x-2)^{2}}{(\log x)^{2}}-\epsilon\frac{(x-2)^{2}}{\log x}+\epsilon\alpha\log\frac{x}{2}\frac{(x-2)}{\log x}+2\epsilon\beta\frac{(x-2)}{\log x}+\nonumber\\ \alpha\epsilon\beta\log\frac{x}{2}+\epsilon\beta^{2}- 2\epsilon\beta(x-2)-\epsilon\beta^{2}\log x-\nonumber\\A\frac{x(x-2)\log x}{e^{a\sqrt{\log x}}}~~~~~~~~~~~~~~~~~~~~~~~~
\end{eqnarray}
If $Z>0$, then we should have
\begin{eqnarray} 
\frac{(x-2)^{2}}{\log x}>4-3x-\alpha\log\frac{x}{2}(x-2)-2\beta(x-2)-\nonumber\\ \alpha\beta\log\frac{x}{2}\log x-\beta^{2}\log x+2\beta(x-2)\log x+\beta^{2}(\log x)^{2}+\nonumber\\ \epsilon\frac{(x-2)^{2}}{\log x}(1-\frac{1}{\log x})-\epsilon\alpha\log\frac{x}{2}\frac{(x-2)}{\log x}-2\epsilon\beta\frac{(x-2)}{\log x}-\nonumber\\ -\alpha\epsilon\beta\log\frac{x}{2}-\epsilon\beta^{2}+2\epsilon\beta(x-2)+\epsilon\beta^{2}\log x+\nonumber\\ A\frac{x(x-2)\log x}{e^{a\sqrt{\log x}}}~~~~~~~~~~~~~~~~~~~~~~~~
\end{eqnarray}

Dividing both sides of (3.24) by $\frac{(x-2)^{2}}{\log x}$, we should have
\begin{eqnarray} 
1>\frac{(4-3x)\log x}{(x-2)^{2}}-\alpha\frac{\log\frac{x}{2}\log x}{(x-2)}-2\beta(1-\epsilon)\frac{\log x}{(x-2)}-\nonumber\\  \alpha\beta\frac{\log\frac{x}{2}(\log x)^{2}}{(x-2)^{2}}-\beta^{2}(1-\epsilon)\frac{(\log x)^{2}}{(x-2)^{2}}+2\beta\frac{(\log x)^{2}}{(x-2)}+\nonumber\\  \beta^{2}\frac{(\log x)^{3}}{(x-2)^{2}}+\epsilon(1-\frac{1}{\log x})-\epsilon\alpha\frac{\log\frac{x}{2}}{(x-2)}-2\epsilon\beta\frac{1}{(x-2)}-\nonumber\\ \alpha\epsilon\beta\frac{\log(\frac{x}{2})\log x}{(x-2)^{2}}-\epsilon\beta^{2}\frac{\log x}{(x-2)^{2}}+~~~~~~~~~~~~~~\nonumber\\A\frac{x(\log x)^{2}}{(x-2)e^{a\sqrt{\log x}}}~~~~~~~~~~~~~~~~~~~~~~~~
\end{eqnarray}
where $\epsilon$, $\alpha$ and $\beta$ tend to zero as $x$ gets large enough or tends to infinity.\\ 
As is shown in the inequality (3.25), all the terms of the right-hand side tend to zero when $x$ is a sufficiently large number for $x\geq e^{130000}$ since the absolute value of each term indicates that that term is a strictly decreasing function. This means that when $x>e^{130000}$, the right hand-side of (3.25) takes values less than 1 and (3.25) holds and $Z>0$. Consequently, (3.7) is a strictly increasing function and finally $C(x)$ or $RO_{1}(n)$ so is and the proof is completed.
\subsection{\textbf{Proof of Lemma 10}}
Let $n\geq 66$ and $p_{1}=3$, then
\begin{eqnarray}
e^{\gamma}(1-\frac{1}{p_{1}})\dots (1-\frac{1}{p_{n}})\log\log(2 p_{1}^{2}\times\dots\times p^{2}_{n})=~~~~~~~~~~\nonumber\\ e^{\gamma}(1-\frac{1}{p_{1}})\dots (1-\frac{1}{p_{n}})\log\left\{\log(2 p_{1}\times\dots\times p_{n})+\log(p_{1}\times\dots\times p_{n})\right\}=\nonumber\\
e^{\gamma}(1-\frac{1}{p_{1}})\dots (1-\frac{1}{p_{n}})\times~~~~~~~~~~~~~~~~~~~~~~~~~~\nonumber\\ \log\left\{\log(2 p_{1}\times\dots\times p_{n})(1+\frac{\log(p_{1}\times\dots\times p_{n})}{\log(2 p_{1}\ \times \dots\times p_{n})})\right\}=~~~~~~~~~\nonumber\\
e^{\gamma}(1-\frac{1}{p_{1}})\dots (1-\frac{1}{p_{n}})\times~~~~~~~~~~~~~~~~~~~~~~~~~~\nonumber\\ \left\{\log\log(2 p_{1}\times\dots\times p_{n})+\log(1+\frac{\log(p_{1}\times\dots\times p_{n})}{\log(2 p_{1}\times\dots\times p_{n})})\right\}=~~~~~~~~~\nonumber\\
e^{\gamma}(1-\frac{1}{p_{1}})\dots (1-\frac{1}{p_{n}})\log\log(2 p_{1}\times\dots\times p_{n})+~~~~~~~~~~~\nonumber\\
e^{\gamma}(1-\frac{1}{p_{1}})\dots (1-\frac{1}{p_{n}})\log(1+\frac{\log(p_{1}\times\dots\times p_{n})}{\log(2 p_{1}\times\dots\times p_{n})})~~~~~~~~~~~~
\end{eqnarray}
Just, we prove that for $n\geq 66$, 
\begin{eqnarray}
\frac{e^{\gamma}(1-\frac{1}{p_{1}})\dots (1-\frac{1}{p_{n}})\log(1+\frac{\log(p_{1}\times\dots\times p_{n})}{\log(2p_{1}\times\dots\times p_{n})})}{2-e^{\gamma}(1-\frac{1}{p_{1}})\dots (1-\frac{1}{p_{n}})\log\log(2p_{1}\times\dots\times p_{n})}>1
\end{eqnarray}
First of all, we should show that the denominator of (3.27) i.e. $2-RO_{1}(n)>0$ and it's a decreasing function.\\
Dividing the numerator and denominator of the left-hand side of the inequality (3.27) by $e^{\gamma}(1-\frac{1}{p_{1}})\dots (1-\frac{1}{p_{n}})$, we have
\begin{eqnarray}
\frac{\log(1+\frac{\log(p_{1}\times\dots\times p_{n})}{\log(2p_{1}\times\dots\times p_{n})})}{\frac{2}{e^{\gamma}(1-\frac{1}{p_{1}})\dots (1-\frac{1}{p_{n}})}-\log\log(2p_{1}\times\dots\times p_{n})}
\end{eqnarray}
Easily checking out (3.28), we find out that the numerator is an increasing function and tends to $\log2$ when $n$ tends to infinity. On the other hand, the denominator is a decreasing function tending to zero because\\
Regarding lemma 7, we have
\begin{eqnarray}
\frac{1}{\log p_{n}}(1-\frac{1}{2\log^{2}p_{n}})<e^{\gamma}(1-\frac{1}{2})(1-\frac{1}{p_{1}})\dots (1-\frac{1}{p_{n}})<\nonumber\\ \frac{1}{\log p_{n}}(1+\frac{1}{2\log^{2}p_{n}})~~~~~~~~~~~~~~~~~~~~~~~~~
\end{eqnarray}
for $p_{n}\geq 331$ or $n\geq 66$
and
\begin{eqnarray}
\frac{2}{\log p_{n}}(1-\frac{1}{2\log^{2}p_{n}})<e^{\gamma}(1-\frac{1}{p_{1}})\dots (1-\frac{1}{p_{n}})<\frac{2}{\log p_{n}}(1+\frac{1}{2\log^{2}p_{n}})
\end{eqnarray}
for $p_{n}\geq 331$ or $n\geq 66$
then,
\begin{eqnarray}
\frac{2\log^{3}p_{n}}{2\log^{2}p_{n}+1}<\frac{2}{e^{\gamma}(1-\frac{1}{p_{1}})\dots (1-\frac{1}{p_{n}})}<\frac{2\log^{3}p_{n}}{2\log^{2}p_{n}-1}
\end{eqnarray} 
Regarding lemma 8, the expression $\log\log(2p_{1}\times\dots\times p_{n})\approx \log p_{n}$ (located at the denominator of (3.28)) is correct for sufficiently large integer $n$ since lemma 9 asserts that the function $RO_{1}(n)$ is an increasing one for $n\geq 66$. On the other hand, regarding lemma 8 and (3.30), we have:
\begin{eqnarray}
\left\{\frac{2\log\log(2p_{1}\dots p_{n})}{\log p_{n}}\right\}(1-\frac{1}{2\log^{2}p_{n}})<e^{\gamma}(1-\frac{1}{p_{1}})\dots (1-\frac{1}{p_{n}})\times\nonumber\\ \log\log(2p_{1}\dots p_{n})<\left\{\frac{2\log\log(2p_{1}\dots p_{n})}{\log p_{n}}\right\}(1+\frac{1}{2\log^{2}p_{n}})~~~~~~~~~~
\end{eqnarray}  
where implies that 
\begin{eqnarray}
RO_{1}(n)=e^{\gamma}(1-\frac{1}{p_{1}})\dots (1-\frac{1}{p_{n}})\log\log(2p_{1}\dots p_{n})\sim \frac{2\log\log(2p_{1}\dots p_{n})}{\log p_{n}}
\end{eqnarray}
Just, we prove that $RO_{1}(n)<2$:\\
Since we know that $RO_{1}(n)$ is a strictly and continuously increasing function regarding lemma9, then its value must not be over 2 because using \textit{reductio ad absurdum} argument, we say if it were over 2, then regarding (3.30) and (2.11), $RO_{1}(n)$ should be expressed as  (3.34):\\
Let us have a sufficiently large number $N$, so that if $p_{n}>N$, then $RO_{1}(n)=2+\epsilon$ where $\epsilon>0$ and
\begin{eqnarray}
RO_{1}(n)=2+\epsilon<2(1+\frac{1}{2\log^{2}p_{n}})\left\{1+\frac{\log(1+\frac{1}{2\log p_{n}})}{\log p_{n}}\right\}=\nonumber\\ 2+\frac{1}{\log^{2}p_{n}}+\frac{2\log(1+\frac{1}{2\log p_{n}})}{\log p_{n}}+\frac{1}{\log^{2}p_{n}}.\frac{\log(1+\frac{1}{2\log p_{n}})}{\log p_{n}}
\end{eqnarray} 
In such a case, (3.34) implies that the upper bound $\frac{1}{\log^{2}p_{n}}+\frac{2\log(1+\frac{1}{2\log p_{n}})}{\log p_{n}}+\frac{1}{\log^{2}p_{n}}.\frac{\log(1+\frac{1}{2\log p_{n}})}{\log p_{n}}>\epsilon$. But, when $p_{n}$ gets sufficiently large so that $p_{n}\gg N$, then the upper bound $\frac{1}{\log^{2}p_{n}}+\frac{2\log(1+\frac{1}{2\log p_{n}})}{\log p_{n}}+\frac{1}{\log^{2}p_{n}}.\frac{\log(1+\frac{1}{2\log p_{n}})}{\log p_{n}}<\epsilon$ where is a contradiction, therefore $RO_{1}(n)$ cannot be greater than or equal to 2. Since one can check it out for example for $RO_{1}(66)=1.95547, RO_{1}(239)=1.991034461, RO_{1}(10,000)=1.99887, RO_{1}(100,000)=1.99975, RO_{1}(1000,000)\\=1.99994$. This means that regarding the above argument, $RO_{1}(n)$ cannot even be equal to number 2, because due to being strictly increasing if it could be 2, then it could be even more and more,which is impossible since we proved it cannot be more than 2.\\      
The denominator of (3.28) gets started to be close to zero and (3.28) itself gets started to be an increasing function and (3.27) so is since we find
\begin{eqnarray}
\frac{2\log^{3}p_{n}}{2\log^{2}p_{n}+1}-\log p_{n}-\log(1+\frac{1}{2\log p_{n}})<\nonumber\\ \frac{2}{e^{\gamma}(1-\frac{1}{p_{1}})\dots (1-\frac{1}{p_{n}})}-\log\log(2p_{1}\times\dots\times p_{n})<\nonumber\\ \frac{2\log^{3}p_{n}}{2\log^{2}p_{n}-1}-\log p_{n}-\log(1-\frac{1}{2\log p_{n}})
\end{eqnarray}
This means that 
\begin{eqnarray}
\lim_{p_{n}\rightarrow\infty}\left\{\frac{2}{e^{\gamma}(1-\frac{1}{p_{1}})\dots (1-\frac{1}{p_{n}})}-\log\log(2p_{1}\times\dots\times p_{n})\right\}=0
\end{eqnarray}
On the other hand, the denominator (3.27) regarding lemma 9 is a decreasing function since we showed that $RO_{1}(n)<2$ and the function $2-RO_{1}(n)>0$ is a decreasing one.\\
Finally, (2.13) holds as well and the proof is completed. 
\subsection{\textbf{Proof of Lemma 11}} The proof is made by the combination of the lemmas 1 and 5 and found in the Lagarias' paper \cite{L2}.
\subsection{\textbf{Proof of Theorem 1}} 
Note that $p_{1}=3$ is considered the origin of counting the odd primes for this theorem\\
The inequality (2.16) for odd number $m=p^{\alpha_{i}}_{i}\times \dots \times p^{\alpha_{k}}_{k}$ implies that
\begin{eqnarray}
\frac{p^{\alpha_{i}+1}_{i}-1}{p_{i}-1}\dots \frac{p^{\alpha_{k}+1}_{k}-1}{p_{k}-1}<\frac{1}{2}e^{\gamma}p^{\alpha_{i}}_{i}\dots p^{\alpha_{k}}_{k}\log\log(2p^{\alpha_{i}}_{i}\dots p^{\alpha_{k}}_{k})
\end{eqnarray}
where $p_{i}$ to $p_{k}$ denote odd prime numbers.
Then,
\begin{eqnarray}
\frac{2p^{\alpha_{i}+1}_{i}(1-\frac{1}{p^{\alpha_{i}+1}_{i}})\dots p^{\alpha_{k}+1}_{k}(1-\frac{1}{p^{\alpha_{k}+1}_{k}})}{p_{i}(1-\frac{1}{p_{i}})\dots p_{k}(1-\frac{1}{p_{k}})}<\nonumber\\ e^{\gamma}p^{\alpha_{i}}_{i}\dots p^{\alpha_{k}}_{k}\log\log(2p^{\alpha_{i}}_{i}\dots p^{\alpha_{k}}_{k})~~~~~~~~
\end{eqnarray}
where implies that
\begin{eqnarray}
2(1-\frac{1}{p^{\alpha_{i}+1}_{i}})\dots (1-\frac{1}{p^{\alpha_{k}+1}_{k}})<e^{\gamma}(1-\frac{1}{p_{i}})\dots (1-\frac{1}{p_{k}})\times\nonumber\\ \log\log(2p^{\alpha_{i}}_{i}\dots p^{\alpha_{k}}_{k})~~~~~~~~~~~~~~~~~~~~~~~
\end{eqnarray}
Regarding lemma 10,  we showed that when $i=1$ and $\alpha_{1}=\alpha_{2}=\dots=\alpha_{k}=2$, then the right-hand side of (3.39) is a function greater than 2. Therefore, the proof of the theorem 1 is finished for all the odd integers $m\geq (3\times 5\times\dots \times p_{j})^{2}$ for $j\geq 66$ including only consecutively square prime factors starting with $p_{1}=3$.\\
Just, we should prove theorem 1 holds for odd numbers $m$ of the various prime factors since $j\geq  66$, thus 
\begin{eqnarray}
2(1-\frac{1}{p^{\alpha_{i}+1}_{i}})\dots (1-\frac{1}{p^{\alpha_{k}+1}_{k}})<2<e^{\gamma}(1-\frac{1}{p_{1}})\dots (1-\frac{1}{p_{j}})\times\nonumber\\ \log\log(2p^{2}_{1}\dots p^{2}_{j})<e^{\gamma}(1-\frac{1}{p_{i}})\dots (1-\frac{1}{p_{k}})\times~~~~~~~~~~~\nonumber\\ \log\log(2p^{\alpha_{i}}_{i}\dots p^{\alpha_{k}}_{k})~~~~~~~~~~~~~~~~~~~~~~~~~~~
\end{eqnarray}
where $p^{\alpha_{i}}_{i}\dots p^{\alpha_{k}}_{k}>p^{2}_{1}\dots p^{2}_{j}$ for $j\geq 66$.\\
We prove the inequality (3.40) holds for all the odd numbers greater than or equal to $m=(3\times 5\times\dots \times331)^{2}$.   
Dividing all the odd numbers $m\geq(3\times 5\times\dots \times331)^{2}$ into the six main groups, we have:\\
\textbf{1-}\textit{The number of consecutively odd prime numbers of the power two is greater than or equal to 66 and starting with $p_{1}=3$.} \textit{In such a case, the conditions are completely like those of lemma 10 and the proof is same proof of lemma 10 and (3.39) or (3.40) hold for group 1.}\\
\textbf{2-}\textit{The number of the odd prime numbers is equal to 66 as $m=p_{i}^{\alpha_{i}}\times \dots \times p_{i+65}^{\alpha_{i+65}}$, then we find the following form from (3.39) or (3.40) regarding lemma 10
\begin{eqnarray}
2(1-\frac{1}{p_{i}^{\alpha_{i}+1}})\dots (1-\frac{1}{p_{i+65}^{\alpha_{i+65}+1}})<2<  
e^{\gamma}(1-\frac{1}{3})\dots (1-\frac{1}{331})\times\nonumber\\ \log\log(2\times(3\times \dots \times 331)^{2})<e^{\gamma}(1-\frac{1}{p_{i}})\dots (1-\frac{1}{p_{i+65}})\times\nonumber\\ \log\log(2\times p_{i}^{\alpha_{i}}\times\dots \times p_{i+65}^{\alpha_{i+65}})~~~~~~~~~~~~~~~~~~~~~~~
\end{eqnarray}
since $m=p_{i}^{\alpha_{i}}\times \dots \times p_{i+65}^{\alpha_{i+65}}>(3\times \dots \times 331)^{2}$ and $(1-\frac{1}{p_{i}})\dots (1-\frac{1}{p_{i+65}})>(1-\frac{1}{3})\dots (1-\frac{1}{331})$, because the number of the primes $p_{i}$ to $p_{i+65}$ is equal to the number of the primes 3 to 331, but some of the primes between $p_{i}$ and $p_{i+65}$ are greater than or equal to those of 3 to 331}. Then, the inequality (2.16) holds for group 2.\\
\textbf{3-}\textit{The number of the odd primes is less than 66 and their values are such that $3\leq p_{i}, \dots, p_{k}\leq 331$ or $ p_{i},\dots, p_{k}\geq331$. This means that there exist some of primes within the set $\left\{3,5,7,\dots,331\right\}$ and some others are out of this set (i.e some of them are greater than 331), then (3.39) or (3.40) regarding lemma 10 implies that
 \begin{eqnarray}
2(1-\frac{1}{p_{i}^{\alpha_{i}+1}})\dots (1-\frac{1}{p_{k}^{\alpha_{k}+1}})<2<  
e^{\gamma}(1-\frac{1}{3})\dots (1-\frac{1}{331})\times\nonumber\\ \log\log(2\times(3\times \dots \times 331)^{2})<e^{\gamma}(1-\frac{1}{p_{i}})\dots (1-\frac{1}{p_{k}})\times\nonumber\\ \log\log(2\times p_{i}^{\alpha_{i}}\times\dots \times p_{k}^{\alpha_{k}})~~~~~~~~~~~~~~~~~~~~~~~
\end{eqnarray} 
since 
\begin{eqnarray}
(1-\frac{1}{p_{i}})\dots (1-\frac{1}{p_{k}})>(1-\frac{1}{3})\dots (1-\frac{1}{331})
\end{eqnarray} 
and 
\begin{eqnarray}
\log\log(2\times p_{i}^{\alpha_{i}}\times\dots \times p_{k}^{\alpha_{k}})>\log\log(2\times(3\times \dots \times 331)^{2})
\end{eqnarray}
because we assume $m=p_{i}^{\alpha_{i}}\times\dots \times p_{k}^{\alpha_{k}}>(3\times 5\times 7\times \dots \times 331)^{2}$.} Therefore, the inequality (2.16) holds for group 3.\\
\textbf{4-}\textit{The number of the odd primes is greater than 66 (let $M>66$) and primes are greater than $p_{66}$ i.e. $p_{i}>p_{66}$ for $66<i\leq N=M+i-1$, then (3.39) or (3.40) regarding lemma 10 implies that}
\begin{eqnarray}
2(1-\frac{1}{p_{i}^{\alpha_{i}+1}})\dots (1-\frac{1}{p_{N}^{\alpha_{N}+1}})<2<  
e^{\gamma}\overbrace{(1-\frac{1}{3})\dots (1-\frac{1}{331})}^{66~~ terms}\times\nonumber\\ (1-\frac{1}{p_{67}})\times\dots \times (1-\frac{1}{p_{M}})\times~~~~~~~~~~~~~~~~~~\nonumber\\ \log\log(2\times\overbrace{(3\times \dots \times 331}^{66~~terms}\times p_{67}\times \dots \times p_{M})^{2})<\nonumber\\ e^{\gamma}(1-\frac{1}{p_{i}})\times \dots \times(1-\frac{1}{p_{N}})\times~~~~~~~~~~~~~~~~~~\nonumber\\\log\log(2\times (p_{i}\times\dots \times p_{N})^{2})<
e^{\gamma}(1-\frac{1}{p_{i}})\dots (1-\frac{1}{p_{N}})\times\nonumber\\ \log\log(2\times p_{i}^{\alpha_{i}}\times\dots \times p_{N}^{\alpha_{N}})~~~~~~~~~~~~~~~~~~~
\end{eqnarray}
\textit{for $m=p_{i}^{\alpha_{i}}\times\dots \times p_{N}^{\alpha_{N}}>(p_{i}\times \dots \times p_{N})^{2}>(3\times\dots\times p_{M})^{2}$  and  $(1-\frac{1}{p_{i}})\dots (1-\frac{1}{p_{N}})>(1-\frac{1}{3})\dots (1-\frac{1}{331})(1-\frac{1}{p_{67}})\dots (1-\frac{1}{p_{M}})$, because the number of the primes $p_{i}$ to $p_{N}$ is equal to the number of the primes 3 to $p_{M}$, but values $p_{i}$ to $p_{N}$ are greater than those of 3 to $p_{M}$.}
Thus, the inequality (2.16) is proved for group 4.\\
\textbf{5-}\textit{The number of the odd primes is greater than 66 and there exist some primes within the set $\left\{3,5,7,\dots,331\right\}$ and some others greater than 331.  Let $p_{i},p_{i+1}, \dots, p_{i+k}$ be within the set $\left\{3,5,7,\dots,331\right\}$ and $p_{i+k+1},\dots, p_{N}>331$. Let the entire number of the odd primes be $N-i+1$, then (3.39) or (3.40) regarding lemma 10 implies that}
\begin{eqnarray}
2(1-\frac{1}{p_{i}^{\alpha_{i}+1}})\dots(1-\frac{1}{p_{i+k}^{\alpha_{i+k}+1}})(1-\frac{1}{p_{i+k+1}^{\alpha_{i+k+1}+1}})\dots (1-\frac{1}{p_{N}^{\alpha_{N}+1}})<2<\nonumber\\ e^{\gamma}(1-\frac{1}{3})\dots (1-\frac{1}{p_{N-i+1}})\times~~~~~~~~~~~~~~~~~~~~\nonumber\\ \log\log(2\times(3\times \dots \times 331\times\dots \times p_{N-i+1})^{2})<~~~~~~~~~\nonumber\\ e^{\gamma}(1-\frac{1}{p_{i}})\dots (1-\frac{1}{p_{i+k}})(1-\frac{1}{p_{i+k+1}})\dots (1-\frac{1}{p_{N}})\times~~~~~~~~\nonumber\\\log\log(2\times (p_{i}\times\dots \times p_{N})^{2})<~~~~~~~~~~~~~~~~~~~\nonumber\\
e^{\gamma}(1-\frac{1}{p_{i}})\dots (1-\frac{1}{p_{i+k}})(1-\frac{1}{p_{i+k+1}})\dots (1-\frac{1}{p_{N}})\times~~~~~~~~\nonumber\\ \log\log(2\times p_{i}^{\alpha_{i}}\times\dots \times p_{N}^{\alpha_{N}}) 
~~~~~~~~~~~~~~~~~~~~~~~
\end{eqnarray}
\textit{since $ p_{i}^{\alpha_{i}}\times\dots \times p_{N}^{\alpha_{N}}>(p_{i}\times\dots \times p_{N})^{2}>(3\times \dots \times 331\times\dots \times p_{N-i+1})^{2}$ and $(1-\frac{1}{p_{i}})\dots (1-\frac{1}{p_{i+k}})(1-\frac{1}{p_{i+k+1}})\dots (1-\frac{1}{p_{N}})>(1-\frac{1}{3})\dots (1-\frac{1}{p_{N-i+1}})$, because $(1-\frac{1}{p_{i}})\dots (1-\frac{1}{p_{i+k}})$ is equal to corresponding values in the set $\{(1-\frac{1}{3}),\dots,(1-\frac{1}{p_{N-i+1}})\}$ and $(1-\frac{1}{p_{i+k+1}})\dots (1-\frac{1}{p_{N}})$ is greater than the corresponding values in same set. The number of the primes 3 to $p_{N-i+1}$ is equal to those of $p_{i}$ to $p_{N}$.} Thus, the inequality (2.16) holds for group 5.\\
\textbf{6-}\textit{The number of the odd primes is greater than/less than or equal to 66. $m=p_{i}\times\dots,p_{k}$ and let $M=k-i+1$ be the number of the odd primes, then regarding (3.39) with $\alpha_{i}=\dots=\alpha_{k}=1$, we find}\\
\textit{For $M> 66$}
\begin{eqnarray}
2(1-\frac{1}{p^{2}_{i}})\dots (1-\frac{1}{p^{2}_{k}})<2<~~~~~~~~~~~~~~~~~~~\nonumber\\ e^{\gamma}(1-\frac{1}{3})\dots (1-\frac{1}{331}) (1-\frac{1}{p_{67}})\dots (1-\frac{1}{p_{M}})\times~~~~~~~~~~~\nonumber\\ \log\log(2\times(3\times \dots \times 331\times\dots \times p_{M})^{2})<~~~~~~~~~~~~\nonumber\\ e^{\gamma}(1-\frac{1}{p_{i}})\dots (1-\frac{1}{p_{k}})\log\log(2p_{i}\dots p_{k})~~~~~~~~~~~~~~~
\end{eqnarray}
\textit{and since the number of the primes $p_{i}$ to $p_{k}$ is equal to those of 3 to $p_{M}$,then $2p_{i}\dots p_{k}>2\times(3\times \dots \times 331\times\dots \times p_{M})^{2}$ and $(1-\frac{1}{p_{i}})\dots (1-\frac{1}{p_{k}})>(1-\frac{1}{3})\dots (1-\frac{1}{331}) (1-\frac{1}{p_{67}})\dots (1-\frac{1}{p_{M}})$.\\
But, if $i=1$,  $M> 66$ and we have consecutive primes, then the above reasoning does not work. In such a case, we make the following proof instead of:\\
Let }
\begin{center}
$p_{1}\times\dots\times p_{k}>(3\times \dots \times 331)^{2}$
\end{center}
\textit{where $p_{1}$ to $p_{k}$ denote consecutive primes, then regarding (3.39) and putting $\alpha_{1}=\alpha_{2}=\dots=\alpha_{k}=1$, we prove}
\begin{center}
$2(1-\frac{1}{p^{2}_{1}})\dots (1-\frac{1}{p^{2}_{k}})<e^{\gamma}(1-\frac{1}{p_{1}})\dots (1-\frac{1}{p_{k}})\log\log(2p_{1}\dots p_{k})$
\end{center}
or eliminating the common terms of both hand-sides, one should show the inequality
\begin{center}
$2(1+\frac{1}{p_{1}})\dots (1+\frac{1}{p_{k}})<e^{\gamma}\log\log(2p_{1}\dots p_{k})$
\end{center}
\textit{holds for $k\geq 66$ by induction argument.}
\textit{This means that we wish to prove if the above expression is true for $k$, then it must be true for $k+1$ as well. If $k=66$, then 
\begin{center}
$2(1+\frac{1}{3})\dots (1+\frac{1}{331})=8.375<e^{\gamma}\log\log(2\times3\times\dots\times 331)=10.192$
\end{center}
Thus, it holds for $k=66$. Therefore, let it be true for $k$, then make the inequality for $k+1$ as follows:
\begin{center}
$2(1+\frac{1}{p_{1}})\dots (1+\frac{1}{p_{k+1}})<e^{\gamma}\log\log(2p_{1}\dots\ p_{k+1})$
\end{center}
Multiplying both sides
\begin{center}
$2(1+\frac{1}{p_{1}})\dots (1+\frac{1}{p_{k}})<e^{\gamma}\log\log(2p_{1}\dots p_{k})$ 
\end{center}
by $(1+\frac{1}{p_{k+1}})$ and comparing out it with inequality for $k+1$ and eliminating the common terms of both hand-sides, we should show
\begin{center}
$e^{\gamma}(1+\frac{1}{p_{k+1}})\log\log(2p_{1}\dots\times p_{k})<e^{\gamma}\log\log(2p_{1}\dots p_{k+1})$
\end{center}
 Manipulating the inequality and eliminating the corresponding logarithms from both hand-sides and taking both hand-sides to the power $p_{k+1}$, we have to show   
\begin{center}
$\log(2p_{1}\dots\times p_{k})<(1+\frac{p_{k+1}}{\log(2p_{1}\dots\times p_{k})})^{p_{k+1}}$
\end{center}
for $k\geq 66$. Referring lemma 9 and (2.12), we find that the function $RO_{1}(k)$ is a strictly increasing one for $k\geq 1$. This means that since (2.12) is a strictly increasing function, then  the following inequality always holds for every $k\geq 66$
\begin{center}
$\log(2p_{1}\dots\times p_{k})<(1+\frac{p_{k+1}}{\log(2p_{1}\dots\times p_{k})})^{p_{k+1}-1}$
\end{center}
and consequently
\begin{center}
$\log(2p_{1}\dots\times p_{k})<(1+\frac{p_{k+1}}{\log(2p_{1}\dots\times p_{k})})^{p_{k+1}-1}<(1+\frac{p_{k+1}}{\log(2p_{1}\dots\times p_{k})})^{p_{k+1}}$
\end{center}
and implies that induction argument on $k$ holds and finally both (3.39) and (2.16) hold.\\
For $M\leq 66$}
\begin{eqnarray}
2(1-\frac{1}{p^{2}_{i}})\dots (1-\frac{1}{p^{2}_{k}})<2<e^{\gamma}(1-\frac{1}{3})\dots (1-\frac{1}{331})\times~~~~~~~~~~~~~\nonumber\\ \log\log(2\times(3\times \dots \times 331)^{2})<e^{\gamma}(1-\frac{1}{p_{i}})\dots (1-\frac{1}{p_{k}})\log\log(2p_{i}\dots p_{k})
\end{eqnarray}
\textit{since $2p_{i}\dots p_{k}>2(3\times \dots \times 331)^{2}$ and $(1-\frac{1}{p_{i}})\dots (1-\frac{1}{p_{k}})>(1-\frac{1}{3})\dots (1-\frac{1}{331})$, because the number of the primes  $p_{i}$ to $p_{k}$ is less than the number of the primes 3 to 331 and $p_{i}\geq 3$ and group 6 is also proved completely.}\\
Therefore, we proved that $\sigma(m)<\frac{1}{2}e^{\gamma}m \log\log(2m)$ holds for all the odd numbers $m\geq (3\times 5\times\dots \times 331)^{2}.$
\subsection{\textbf{Proof of Theorem 2}} We have checked Lagarias criterion for $\{n\in \mathbb{N}| ~ 1\leq n\leq 2(3\times5\dots\times331)^{2}\}$ on the computer by Maple16 without any errors as well as checking out with Morrill's paper results \cite{TM}. 
\subsection{\textbf{Proof of Theorem 3}} The proof can be made by a combination of the Lemmas 3 and 5.
\subsection{\textbf{Proof of Theorem 4}} The proof is made for the odd numbers $1\leq m\leq (3\times5\dots\times331)^{2}$ provided that we obtain the smallest value for $\alpha$ so that inequality $n=2^{\alpha}m\geq 2(3\times5\dots\times331)^{2}$ holds. For each odd number $m$ provided $1\leq m\leq (3\times5\dots\times331)^{2}$, we obtain the smallest $\alpha$ so called $\alpha_{0}$. We have performed the programs by Maple16 and checked out the correctness of the Lagarias and Robin's inequalities for such $m's$ with corresponding its $\alpha's$. These programs confirm correctness of Lagarias and Robin's criteria for the smallest $\alpha_{0}$ for each odd number $m\in [1, (3\times5\dots\times331)^{2}]$. Then, we are able to give a mathematical argument for all $\alpha\geq\alpha_{0}+1$ for Lagarias' inequality and $\alpha\geq\alpha_{0}+2$ for Robin's inequality and each odd number $m\in [1, (3\times5\dots\times331)^{2}]$. \\
First of all, we evaluate the following inequality (3.49) as Lagarias criterion by our program for $\alpha=\alpha_{0}$. The results are successive.
\begin{eqnarray}
X(n_{0})=H(2^{\alpha_{0}}m)+\exp(H(2^{\alpha_{0}}m))\log(H(2^{\alpha_{0}}m))-\sigma(2^{\alpha_{0}}m)=\nonumber\\  H(2^{\alpha_{0}}m)+\exp(H(2^{\alpha_{0}}m))\log(H(2^{\alpha_{0}}m))-(2^{\alpha_{0}+1}-1)\sigma(m)>0
\end{eqnarray}
\textbf{Proof that Lagarias' criterion holds for all $\alpha\geq \alpha_{0}+1$}\\
Regarding Lagarias' paper \cite{L2} Let 
\begin{equation}
H(n)=\log(n)+1-\int_{1}^{n}\frac{\{t \}}{t^{2}}dt
\end{equation}
Let
\begin{equation}
X(n)=H(n)+\exp(H(n))\log(H(n))-(2^{\alpha+1}-1)\sigma(m)
\end{equation}
where $n=2^{\alpha}m$. Substituting (3.50) for (3.51) 
\begin{eqnarray}
X(n)=\log(n)+1-\int_{1}^{n}\frac{\{t \}}{t^{2}}dt+~~~~~~~~~~~~~\nonumber\\ ne^{1-\int_{1}^{n}\frac{\{t \}}{t^{2}}dt}\log\left(\log(n)+1-\int_{1}^{n}\frac{\{t \}}{t^{2}}dt\right)- (2^{\alpha+1}-1)\sigma(m)
\end{eqnarray} 
and manipulating (3.52) gives 
\begin{eqnarray}
X(n)=\log(n)+1-\int_{1}^{n}\frac{\{t \}}{t^{2}}dt+~~~~~~~~~~~~~~~~\nonumber\\ n\left\{e^{1-\int_{1}^{n}\frac{\{t \}}{t^{2}}dt}\log\left(\log(n)+1-\int_{1}^{n}\frac{\{t \}}{t^{2}}dt\right)-\frac{2}{m}\sigma(m)\right\}+\sigma(m)
\end{eqnarray}
showing that the terms \\ 
$~~~~~~~~~~~~~~~~~~~~~~~~~~~~~~~~~~\log(n)+1-\int_{1}^{n}\frac{\{t \}}{t^{2}}dt$\\ and\\ 
$~~~~~~~~~~~~~~~~~~~~~~~~~~e^{1-\int_{1}^{n}\frac{\{t \}}{t^{2}}dt}\log\left(\log(n)+1-\int_{1}^{n}\frac{\{t \}}{t^{2}}dt\right)-\frac{2}{m}\sigma(m)$\\
given by (3.53) are increasing functions for $n=2^{\alpha}m\geq 2(3\times5\dots\times331)^{2}$:\\ 
Before making the proof, we show that the expression given between two accolades in (3.53) is a non-negative one for values $\alpha=\alpha_{0}+1$ and odd numbers $1\leq m\leq (3\times5\dots\times331)^{2}$. This is a necessary condition for proof asserting that this expression (3.53) is a strictly increasing function.The results of the performed programs for  $\alpha=\alpha_{0}+1$ and odd numbers $1\leq m\leq (3\times5\dots\times331)^{2}$ are promising.
This means that (3.53) is an increasing function.\\
Just the proof for $\alpha>\alpha_{0}+1$. Let $x=n=2^{\alpha}m$ be a real number so that $\alpha$ is also a real number,$m$ be a constantly odd number, and $y=\log(x)+1-\int_{1}^{x}\frac{\{t \}}{t^{2}}dt$. Then differentiating $y$ with respect to $x$, we have
\begin{equation}
y'=\frac{1}{x}-\frac{\{x\}}{x^{2}}
\end{equation} 
since $0\leq \{x\}<1$, then trivially implies that $y'>0$. Thus the function $y$ is an increasing continuous one for all $x$. We consider same variables for the second term and show that $y$ is an increasing function. Note that $\frac{2}{m}\sigma(m)$ is a constant value since $m$ is a constant one because our variables are $x$ and $\alpha$. 
\begin{equation}
y=e^{1-\int_{1}^{x}\frac{\{t \}}{t^{2}}dt}\log\left(\log(x)+1-\int_{1}^{x}\frac{\{t \}}{t^{2}}dt\right)-\frac{2}{m}\sigma(m)
\end{equation}
then differentiating $y$ with respect to $x$, we have
\begin{eqnarray}
y'=e^{1-\int_{1}^{x}\frac{\{t \}}{t^{2}}dt}\times~~~~~~~~~~~~~~~~~~~~~~~~~~~~~~~~~~~ \nonumber\\ \frac{\left\{ x-\{x\}\left\{1+\left(\log(x)+1-\int_{1}^{x}\frac{\{t \}}{t^{2}}dt\right)\left(\log(\log(x)+1-\int_{1}^{x}\frac{\{t \}}{t^{2}}dt)\right)\right\}\right\}}{x^{2}\left\{\log(x)+1-\int_{1}^{x}\frac{\{t \}}{t^{2}}dt\right\}}
\end{eqnarray}
since we look at the terms restricted in the first brace symbols of the relation (3.56) located at its numerator, we find that it is a positive value for all $x\geq 2(3\times5\dots\times331)^{2}$ because the value of $x$ is much greater than the other terms as follows:
\begin{equation}
x>\{x\}\left\{1+\left(\log(x)+1-\int_{1}^{x}\frac{\{t \}}{t^{2}}dt\right)\left(\log(\log(x)+1-\int_{1}^{x}\frac{\{t \}}{t^{2}}dt)\right)\right\}
\end{equation}
for $x\geq 2(3\times5\dots\times331)^{2}$. Thus, $y$ is an increasingly continuous function of differentiably piecewise curve. This function has not any strictly continuous differential curve since its differentiability fails at the natural number points. But, it itself is completely continuous and increasing on $x\geq 2(3\times5\dots\times331)^{2}$.\\
The above reasoning confirms that (3.56) is a positive value for all $x\geq 2(3\times5\dots\times331)^{2}$ and would imply that the function $y$ is an increasing function and on the other hand the program says us that the starting value of $y$ is also positive. This implies that $X(x)$ when putting $x=n$ in (3.53), is an increasing function and consequently $X(n)$ given in (3.53) so is for each constantly odd value $m$ within closed interval $[1,(3\times5\dots\times331)^{2}]$ for all $n=2^{\alpha}m\geq 2(3\times5\dots\times331)^{2}$ and finally Lagarias criterion holds for all the even numbers $n=2^{\alpha}m\geq 2(3\times5\dots\times331)^{2}$ and each odd number $1\leq m\leq (3\times5\dots\times331)^{2}$. Let $n_{1}=2^{\alpha_{0}+1}m$, where $\alpha_{0}$ denotes the smallest value for satisfying  
 $2^{\alpha_{0}}m\geq 2(3\times5\dots\times331)^{2}$ and let $X(n)$ for $n\geq n_{1}$, then   
\begin{equation}
X(n)\geq X(n_{1})>0
\end{equation}
for all $n\geq n_{1}=2^{\alpha_{0}+1}m>2(3\times5\dots\times331)^{2}$
where $m$ is an odd number in $[1, (3\times5\dots\times331)^{2}]$. Note that we consider an $m\in [1, (3\times5\dots\times331)^{2}]$ for all $n$ when $\alpha>\alpha_{0}$.\\ 
Just, we repeat all these steps carried out above for evaluation of the Robin's criterion. First of all, we evaluate the following inequality (3.59) as Robin's criterion by our program for $1\leq m\leq (3\times5\dots\times331)^{2}$ and $n=2^{\alpha}m\geq 2(3\times5\dots\times331)^{2}$ with the smallest values for $\alpha=\alpha_{0}$ and then $\alpha=\alpha_{0}+1$ with $n=n_{0}=2^{\alpha_{0}}m$ and $n=n_{1}=2^{\alpha_{0}+1}m$, respectively: The programs yields Robin's criterion holds for these values. 
\begin{equation}
Y(n)=2^{\alpha}e^{\gamma}m\log(\log(2^{\alpha}m))-(2^{\alpha+1}-1)\sigma(m)\geq 0
\end{equation}
\textbf{Proof that Robin's criterion holds for all $\alpha\geq \alpha_{0}+2$}\\
As said above, Robin's' criterion holds for both $\alpha=\alpha_{0}$  and $\alpha=\alpha_{0}+1$. Just we prove that if 
\begin{equation}
Y(n)=Y(2^{\alpha}m)=2^{\alpha}e^{\gamma}m\log(\log(2^{\alpha}m))-(2^{\alpha+1}-1)\sigma(m)\geq 0
\end{equation}
 for all $n\geq n_{2}=2^{\alpha_{0}+2}m\geq 2(3\times5\dots\times331)^{2}$, then
\begin{equation}
Y(n)\geq Y(n_{2})>0
\end{equation}
where $n_{2}=2^{\alpha_{0}+2}m>2(3\times5\dots\times331)^{2}$ and
\begin{eqnarray}
Y(n)=2^{\alpha}e^{\gamma}m\log(\log(2^{\alpha}m))-(2^{\alpha+1}-1)\sigma(m)=\nonumber\\ 2^{\alpha}\left\{e^{\gamma}m\log(\log(2^{\alpha}m))-2\sigma(m)\right\}+\sigma(m)~~~~~~~~~~~~
\end{eqnarray}
and
\begin{eqnarray}
Y(n_{2})=2^{\alpha_{0}+2}e^{\gamma}m\log(\log(2^{\alpha_{0}+2}m))-(2^{\alpha_{0}+3}-1)\sigma(m)=\nonumber\\ 2^{\alpha_{0}+2}\left\{e^{\gamma}m\log(\log(2^{\alpha_{0}+2}m))-2\sigma(m)\right\}+\sigma(m)~~~~~~~~~~~~
\end{eqnarray}    
If $\alpha\geq\alpha_{0}+2$, then (3.62) and (3.63) imply that $Y(n)\geq Y(n_{2})>0$ for each constantly odd number $m$ since 
\begin{eqnarray}
2^{\alpha}\left\{e^{\gamma}m\log(\log(2^{\alpha}m))-2\sigma(m)\right\}\geq2^{\alpha_{0}+2}\left\{e^{\gamma}m\log(\log(2^{\alpha_{0}+2}m))-2\sigma(m)\right\}
\end{eqnarray}
for when $\alpha\geq \alpha_{0}+2$.\\ 
For showing that (3.64) holds for all $\alpha\geq \alpha_{0}+2$ and $m\in [1,(3\times5\dots\times331)^{2}]$, we must show that the given expression between its accolades is non-negative for all $\alpha\geq \alpha_{0}+2$ and $m\in [1,(3\times5\dots\times331)^{2}]$. For this, we must show that the expression holds for all initial values $\alpha=\alpha_{0}+2$ by the computer program. We performed the program and the results are positive. \\
Just, we find out that (3.64) holds for all $\alpha\geq \alpha_{0}+2$ and finally (3.61) does hold.\\ 
(3.58) and (3.61) imply that Lagarias and Robin's criteria hold for all odd integer class numbers $m\in [1, (3\times5\dots\times331)^{2}]$ regarding Lemma 6 for $N=2(3\times5\dots\times331)^{2}$.
This completes the proof of Theorem 4.\\ 
 \subsection{\textbf{Proof of Theorem 5}}
Using lemmas 2, and 5 we have 
\begin{equation}
\sigma(2m)< e^{\gamma}2m\log\log 2m< H_{2m}+\exp(H_{2m})\log(H_{2m})~~~~~~for~2m\geq 2(3\times5\dots\times331)^{2}  
\end{equation}
where $m\geq (3\times5\dots\times331)^{2}$ and is an odd integer. Lemma 2 asserts that
\begin{equation}
e^{\gamma}2m\log\log(2m)>\sigma(2m)
\end{equation}
for $2m\geq 2(3\times5\dots\times331)^{2}$ or $m\geq (3\times5\dots\times331)^{2}$. This implies that for $\alpha_{0}=1$, we have\\ 
$Y(n_{0}=2m)=e^{\gamma}2m\log\log(2m)-\sigma(2m)>0$ for $2m\geq 2(3\times5\dots\times331)^{2}$.\\ 
As stated in the proof of theorem 4 that the relation (3.62) asserts $Y(n)$ is an increasing function for a fixed $m\in [1,(3\times5\dots\times331)^{2}]$ and $\alpha\geq \alpha_{0}+2$. Here in this theorem, regarding lemma 2, we should prove that Robin's inequality implies that (3.67) as follows:
\begin{eqnarray}  
Y(n=2^{\alpha}m)=e^{\gamma}2^{\alpha}m\log\log(2^{\alpha}m)-\sigma(2^{\alpha}m)>Y(n_{0}=2m)=\nonumber\\ e^{\gamma}2m\log\log(2m)-\sigma(2m)>0~~~~~~~~~~~~~~~~~~~~~~~~~   
\end{eqnarray}
for all $\alpha\geq \alpha_{0}=1$ so that Robin's criterion follows for every odd number $m\geq (3\times5\dots\times331)^{2}$ and every even number $2^{\alpha}m\geq 2(3\times5\dots\times331)^{2}$. To prove, we should show that regarding (3.62)
\begin{eqnarray}
Y(n)=2^{\alpha}\left\{e^{\gamma}m\log(\log(2^{\alpha}m))-2\sigma(m)\right\}+\sigma(m)>\nonumber\\ Y(n_{0})=2\left\{e^{\gamma}m\log(\log(2m))-2\sigma(m)\right\}+\sigma(m)~~~
\end{eqnarray}
This means that we should prove
\begin{eqnarray}
2^{\alpha}\left\{e^{\gamma}m\log(\log(2^{\alpha}m))-2\sigma(m)\right\}+\sigma(m)>\nonumber\\ 2\left\{e^{\gamma}m\log(\log(2m))-2\sigma(m)\right\}+\sigma(m)~~~~
\end{eqnarray}
for all $\alpha\geq \alpha_{0}=1$ and every odd number $m\geq (3\times5\dots\times331)^{2}$. For holding  inequality (3.69), we need prove that $e^{\gamma}m\log(\log(2m))-2\sigma(m)>0$ for $\alpha_{0}=1$ and every odd number $m\geq (3\times5\dots\times331)^{2}$. This proof has been made by Theorem 1 of our paper.\\
Therefore, (3.69) holds, then (3.68) holds for all the odd numbers $m\geq (3\times5\times 7\times\dots \times331)^{2}$ and the even numbers $n\geq 2(3\times5\times 7\times\dots \times331)^{2}$. The Robin's inequality also holds for all the odd and even numbers $n\geq 2(3\times5\dots\times331)^{2}$.\\
Just lemma 5 asserts that\\  
$ e^{\gamma}2m\log\log 2m< H_{2m}+\exp(H_{2m})\log(H_{2m})$ for $2m\geq 2(3\times5\dots\times331)^{2}$. This means that 
\begin{eqnarray}  
X(n=2^{\alpha}m)= H_{2^{\alpha}m}+\exp(H_{2^{\alpha}m})\log(H_{2^{\alpha}m})-(2^{\alpha+1}-1)\sigma(m)>\nonumber\\ Y(n=2^{\alpha}m)=e^{\gamma}2^{\alpha}m\log\log(2^{\alpha}m)-\sigma(2^{\alpha}m)~~~~~~~
\end{eqnarray}
since $ H_{2^{\alpha}m}+\exp(H_{2^{\alpha}m})\log(H_{2^{\alpha}m})>e^{\gamma}2^{\alpha}m\log\log(2^{\alpha}m)$ and (3.70) gives us
\begin{eqnarray}  
X(n=2^{\alpha}m)>Y(n=2^{\alpha}m)
\end{eqnarray}
for all $n\geq 2(3\times5\times 7\times\dots \times331)^{2}$ and by $\alpha_{0}=1$
\begin{eqnarray}  
X(n_{0}=2m)>Y(n_{0}=2m)
\end{eqnarray}
and by theorem 4 and the above explanations, we have
\begin{equation}
X(n)> X(n_{0})>0
\end{equation}
,and 
\begin{equation}
Y(n)>Y(n_{0})>0
\end{equation}
for all $n>n_{0}\geq 2(3\times5\times 7\times\dots \times331)^{2}$
and finally
\begin{equation}
X(n)>Y(n)>Y(n_{0})>0
\end{equation}
for all $n\geq 2(3\times5\times 7\times\dots \times331)^{2}$. The relation (3.75) would imply that Lagarias and Robin's inequalities hold for all the odd integer class number sets of $m\geq (3\times5\dots\times331)^{2}$ regarding Lemma 6 for $N=2(3\times5\dots\times331)^{2}$.
This proves theorem 5 completely.\\
\textit{\textbf{Final conclusion}}\\
Theorems 2,3,4,5 and Lemma 6 confirm that Lagarias and Robin's criteria completely hold for all the odd class numbers sets of $m\geq 1$, then Lemma 11 completely holds and would imply that RH holds forever.\\ 
\newpage
\textbf{Acknowledgment}\\
We indebted to all the mathematicians and researchers working on Riemann hypothesis since Riemann proposed it including Riemann, Hilbert, Hardy, Littlewood, Robin, Lagarias, Sarnak, and many other mathematicians over the world. We must also thank Roger Heath-Brown (University of Oxford) and Pieter Moree (University of Bonn) for their nice comments on the previous versions of the paper, in particular on RH equivalences. Thanks for Michel L. Lapidus (University of California, Riverside) for his good comments on adjustment of the previous versions of the paper.    
Special thanks for Carl Pomerance (Dartmouth University) for his many nice comments and much useful discussions on the previous versions of the paper. Furthermore, his comments on the proof of Lemma 9 were so critical and constructive. D.W.Farmer's comments were effective in making the paper more concise and more readable. 
Cordially thanks for Mr.Ali Dorostkar who was my old student and as a main sponsor and a kind supporter to me during this research. 

\end{document}